\documentstyle [11pt,twoside]   {article}
  \newtheorem{theorem}{Theorem}
  \newtheorem{lemma}{Lemma}

  \newtheorem{remark}{Remark}

  \newcommand {\pf}  {\mbox{\sc Proof. \,\,}}

  \newcommand {\qed} {\null \hfill \rule{2mm}{2mm}}

\def  \Z {{\bf Z}}

\def  \Q {{\bf Q}}
\def  \K {{\bf K}}

\begin {document}

\title{{\Large{\bf On the size of Diophantine $m$-tuples }}}

\author
{\normalsize ANDREJ DUJELLA }

\date{}
\maketitle

\begin{abstract}
Let $n$ be a nonzero integer and assume that a set $S$ of 
positive integers has the property that $xy+n$ is a perfect 
square whenever $x$ and $y$ are distinct elements of $S$. 
In this paper we find some upper bounds for the size of the 
set $S$. We prove that if $|n|\leq 400$ then $|S|\leq 32$, 
and if $|n|>400$ then $|S|<267.81 \,\log{|n|} \,(\log\log{|n|})^2$. 
The question whether there exists an absolute bound 
(independent on $n$) for $|S|$ still remains open. 
\end{abstract}

\footnotetext{
{\it 2000 Mathematics Subject Classification.} Primary 11D09, 11D45; 
Secondary 11J68, 11N36.

{\it Key words.} Diophantine $m$-tuple, simultaneous approximation of 
algebraic numbers, gap principle, large sieve}

\section{Introduction}
Let $n$ be a nonzero integer. A set of $m$ positive integers
$\{a_1,a_2,\ldots,a_m\}$ is said to have the property $D(n)$ if
$a_ia_j+n$ is a perfect square for all $1\leq i<j\leq m$.
Such a set is called a Diophantine $m$-tuple (with the property
$D(n)$), or $P_n$-set of size $m$.

Diophantus found the quadruple $\{1,\,33,\,68,\,105\}$ with the property
\linebreak $D(256)$. The first Diophantine quadruple with the property $D(1)$, 
the set $\{1,\,3,\,8,\,120\}$, was found by Fermat (see \cite{Dic,Dio}). 
Baker and Davenport \cite{B-D} proved
that this Fermat's set cannot be extended to the Diophantine quintuple,
and a famous conjecture is that there does not exist a Diophantine
quintuple with the property $D(1)$. The theorem of Baker and Davenport
has been recently generalized to several parametric families
of quadruples \cite{D-pdeb1,D-pams,DP1}, but the conjecture
is still unproved.

On the other hand, there are examples of Diophantine quintuples and
sextuples like $\{1,\,33,\,105,\,320,\,18240\}$ with the property $D(256)$
\cite{D-acta2} and \linebreak
$\{99,\,315,\,9920,\,32768,\,44460,\,19534284\}$ with
the property $D(2985984)$ \cite{Gibbs1}.

The purpose of this paper is to find some upper bounds for the
numbers $M_n$ defined by
\[ M_n=\sup \{ |S| \,:\, \mbox{$S$ has the property $D(n)$} \}, \]
where $|S|$ denotes the number of elements in the set $S$.

Considering congruences modulo $4$, it is easy to prove that
$M_{4k+2}=3$ for all integers $k$ (see \cite{Bro,G-S,M-R}). 
In \cite{D-acta1} we proved that if $n\not\equiv 2\!\!\pmod{4}$ and 
$n\not\in\{-4,\,-3,\,-1,\,3,\,5,\,8,\,12,\,20\}$, then 
$M_n\geq 4$. 
Recently, we were able to prove that $M_1\leq 8$ (see \cite{D-9}).
(As we said before, the conjecture is that $M_1=4$.)
Since a set with the property $D(4)$ may contain at most two odd elements,
this result implies $M_4\leq 10$.

Since the number of integer points on the elliptic curve
\begin{equation} \label{1}
y^2=(a_1x+n)(a_2x+n)(a_3x+n)
\end{equation}
is finite, we conclude that there does not exist an infinite set
with the property $D(n)$. However, bounds for the size \cite{Bak} and
for the number \cite{Sch} of solutions of (\ref{1}) depend not only
on $n$ but also on $a_1,a_2,a_3$.

On the other hand, we may consider the hyperelliptic curve
\begin{equation} \label{2}
y^2=(a_1x+n)(a_2x+n)(a_3x+n)(a_4x+n)(a_5x+n)
\end{equation}
of genus $g=2$. Caporaso, Harris and Mazur \cite{C-H-M} proved that
the Lang conjecture on varieties of general type implies that for $g\geq 2$
the number $B(g,\K)=\max_{C} |C(\K)|$ is finite. Here $C$ runs over all
curves of genus $g$ over a number field $\K$, and $C(\K)$ denotes the set of
all $\K$-rational points on $C$.
However, even the question whether $B(2,\Q)<\infty$ is still open.
An example of Keller and Kulesz \cite{K-K} shows that $B(2,\Q)\geq 588$
(see also \cite{Elk,Sta}).
Since $M_n\leq 5+B(2,\Q)$\, (by \cite{P-H-Z} we have also 
$M_n\leq 4+B(4,\Q(\sqrt{n})$), we see that the
Lang conjecture implies that
\[ M=\sup \{M_n \,:\, n\in \Z\setminus \{0\} \} \]
is finite.

At present we are able to prove only the weaker result that
$M_n$ is finite for all $n\in \Z\setminus  \{0\}$. In the proof of this
result we will try to estimate the number of "large" (greater than $|n|^3$),
"small" (between $n^2$ and $|n|^3$) and "very small" (less that $n^2$)
elements of a set with the property $D(n)$. Let us introduce the
following notation:
\begin{eqnarray*}
A_n \!\!&=&\!\! \sup \{ |S\cap [|n|^3, +\infty \rangle |\,:\,
\mbox{$S$ has the property $D(n)$} \}, \\
B_n \!\!&=&\!\! \sup \{ |S\cap \langle n^2, |n|^3 \rangle |\,:\,
\mbox{$S$ has the property $D(n)$} \}, \\
C_n \!\!&=&\!\! \sup \{ |S\cap [1,n^2] |\,:\,
\mbox{$S$ has the property $D(n)$} \}.
\end{eqnarray*}

In estimating the number of "large" elements, we used a theorem of
Bennett \cite{Ben} on simultaneous approximations of algebraic numbers
and a very useful gap principle. We proved

\begin{theorem} \label{tm:1}
$\,\,\,A_n\leq 21$\, for all nonzero integers $n$.
\end{theorem}

For the estimate of the number of "small" elements we used a
"weak" variant of the gap principle and we proved

\begin{theorem} \label{tm:2}
$\,\,\,B_n < 0.65 \,\log{|n|} +2.24$\, for all nonzero integers $n$.
\end{theorem}

Finally, in the estimate of the number of 
"very small" elements we used a large sieve
method due to Gallagher \cite{Gal} and we proved

\begin{theorem} \label{tm:3}
$\,\,\,C_n < 265.55 \,\log{|n|}\,(\log\log{|n|})^2 +9.01 \,\log\log{|n|}$\,
for $|n|\geq 400$.
\end{theorem}

Since we checked that $C_n\leq 5$ for $|n|\leq 400$, we may combine
Theorems \ref{tm:1}, \ref{tm:2} and \ref{tm:3} to obtain

\begin{theorem} \label{tm:4}
\begin{eqnarray*}
M_n \!\!&\leq&\!\! 32 \quad \mbox{for $|n|\leq 400$}, \\
M_n \!\!&<&\!\! 267.81 \,\log{|n|} \,(\log\log{|n|})^2 \quad
\mbox{for $|n|\geq 400$}.
\end{eqnarray*}
\end{theorem}

\section{Large elements}
Assume that the set $\{a,b,c,d\}$ has the property $D(n)$.
Let $ab+n=r^2$, $ac+n=s^2$, $bc+n=t^2$, where $r,s,t$ are nonegative
integers. Eliminating $d$ from the system
\[ ad+1=x^2, \quad bd+1=y^2, \quad cd+1=z^2 \]
we obtain the following system of Pellian equations
\begin{eqnarray}
az^2-cx^2 \!\!&=&\!\! n(a-c), \label{3} \\
bz^2-cy^2 \!\!&=&\!\! n(b-c). \label{4}
\end{eqnarray}

We will apply the following theorem of Bennett \cite{Ben} on
simultaneous approximations of square roots of two rationals which
are very close to $1$.

\begin{theorem}[\cite{Ben}]  \label{tm:B}
If $c_i$, $p_i$, $q$ and $L$ are integers for $0\leq i \leq 2$,
with $c_0 < c_1 < c_2$, $c_j=0$ for some $0\leq j\leq 2$, $q$ nonzero
and $L>M^9$, where
\[ M= \max \{|c_0|,\,|c_1|,\,|c_2|\},  \]
then we have
\[ \max_{0\leq i \leq 2}
    \Big\{ \Big| \sqrt{1+\frac{c_i}{L}} - \frac{p_i}{q} \Big| \Big\}
    > (130 L \gamma)^{-1} q^{-\lambda} \]
where
\[ \lambda= 1+ \frac{\log (33L \gamma)}
   {\log \Big( 1.7 L^2 \prod_{0\leq i<j\leq 2} (c_i-c_j)^{-2} \Big)} \]
and
\[ \gamma = \left\{ \begin{array}{ll}
   \frac{(c_2-c_0)^2(c_2-c_1)^2} {2c_2-c_0-c_1} &
   \mbox{if \,\,$c_2-c_1 \geq c_1 -c_0$}, \\
    \frac{(c_2-c_0)^2(c_1-c_0)^2} {c_1+c_2-2c_0} &
    \mbox{if \,\,$c_2-c_1 < c_1 -c_0$}.
   \end{array}
   \right. \]
\end{theorem}

We will apply Theorem \ref{tm:B} to the numbers
\begin{eqnarray*}
{\theta}_1 \!\!&=&\!\! \frac{s}{a} \sqrt{ \frac{a}{c}} =
\sqrt{ \frac{ac+n}{ac}}=  \sqrt{1+ \frac{n}{ac}}
=\sqrt{1+\frac{nb}{abc}} \,, \\
{\theta}_2 \!\!&=&\!\! \frac{t}{b} \sqrt{ \frac{b}{c}}
   =\sqrt{ \frac{bc+n}{bc}}=  \sqrt{1+ \frac{n}{bc}} =
 \sqrt{1+ \frac{na}{abc}} \,.
\end{eqnarray*}

\begin{lemma} \label{l:1}
Assume that $a<b<c$ and $ac>n$. Then all positive integer solutions
$x,y,z$ of the system (\ref{3}) and (\ref{4}) satisfy
\[ \max \Big( |{\theta}_1 -\frac{sbx}{abz}|,\,
|{\theta}_2-\frac{zay}{abz}| \Big) < \frac{c\cdot |n|}{a} z^{-2} . \]
\end{lemma}

\pf
We have
\[ \Big| \frac{s}{a}\sqrt{\frac{a}{c}} -\frac{sbx}{abz}\Big| =
\frac{s}{az\sqrt{c}} |z\sqrt{a}-x\sqrt{c}|=
\frac{s}{az\sqrt{c}} \cdot \frac{|n(c-a)|}{z\sqrt{a}+x\sqrt{c}} \,.\]
If $n<0$, then $s=\sqrt{ac-|n|}<\sqrt{ac}$ and we obtain
\[ |{\theta}_1 -\frac{sbx}{abz}| < \frac{\sqrt{ac}\cdot |n|\cdot c}
{a\sqrt{ac} z^2}= \frac{c|n|}{a}z^{-2}. \]
If $n>0$, then $x\sqrt{c}>z\sqrt{a}$ and we obtain
\[ |{\theta}_1 -\frac{sbx}{abz}|< \frac{\sqrt{ac+n}\cdot n \cdot c}
{2a\sqrt{ac}z^2} =\sqrt{1+\frac{n}{ac}} \cdot \frac{cn}{2a}z^{-2}  <
\frac{cn}{a}z^{-2}. \]

In the same manner, we obtain $|{\theta}_2-\frac{tay}{abz}|<
\frac{c|n|}{b}z^{-2} <\frac{c|n|}{a}z^{-2}$.
\qed

\begin{lemma} \label{l:2}
Let $\{a,b,c,d\}$, $a<b<c<d$, be a Diophantine quadruple with the
property $D(n)$. If $c>b^{11}|n|^{11}$, then $d\leq c^{131}$.
\end{lemma}

\pf
Let $r,s,t,x,y,z$ be defined as in the beginning of this section.
We will apply Theorem \ref{tm:B} with $\{c_0,c_1,c_2\} =\{0, na,nb\}$,
$L=abc$, $M=|nb|$, $q=abz$, $p_1=sbx$, $p_2=tay$.
Since $abc>|n|^9 b^9$, the condition $L>M^9$ is satisfied. For the
quantity $\gamma$ from Theorem \ref{tm:B} we have
\,$\gamma=\frac{b^2(b-a)^2}{2b-a} |n|^3$\, if $b\geq 2a$ and
\,$\gamma=\frac{a^2b^2}{a+b} |n|^3$\, if $a<b\leq 2a$. In both cases we
have
\[ \frac{b^3}{6} |n|^3 \leq \gamma < \frac{b^3}{2} |n|^3. \]
For the quantity $\lambda$ from Theorem \ref{tm:B} we have
\[ \lambda = 1+ \frac{\log (33abc\gamma)}
  {\log (1.7 c^2(b-a)^{-2} n^{-6})} = 2 - {\lambda}_1 \,,\]
where
\[ {\lambda}_1 = \frac{\log \frac{1.7c}{33ab(b-a)^{2} n^6\gamma}}
   {\log (1.7 c^2(b-a)^{-2} n^{-6})} \,.\]
Theorem \ref{tm:B} and Lemma \ref{l:1} imply
\[ \frac{c|n|}{az^2} > (130abc\gamma)^{-1} (abz)^{{\lambda}_1 -2}
> (130abc\gamma)^{-1} a^{-2} b^{-2} z^{{\lambda}_1-2} .\]
This implies
\[ z^{{\lambda}_1} < 130 a^2b^3c^2 |n| \gamma  \]
and
\begin{equation}  \label{5}
\log z < \frac{\log{(130 a^2b^3c^2 |n| \gamma)}
 \log{(1.7c^2(b-a)^{-2} n^{-6})} }
  {\log (\frac{1.7c}{33ab(b-a)^{2} n^6 \gamma}) } \,.
\end{equation}

Let us estimate the right hand side of (\ref{5}). We have
\[ 130 a^2b^3c^2 |n|\gamma < 65 a^2b^6c^2 n^4 < c^3 \cdot
\frac{65a^2}{b^5 |n|^7} < c^3, \]
unless $n=-1$, $a=1$, $b=2$. However, in \cite{D-ost} it was proved
that the Diophantine pair $\{1,2\}$ with the property $D(-1)$
cannot be extended to a Diophantine quadruple.

The same result implies also that if $|n|=1$, then $b-a>1$. Therefore
\[ 1.7 c^2 (b-a)^{-2} n^{-6} <  c^2 .\]

Finally,
\[ \frac{1.7c}{33ab(b-a)^2 n^6 \gamma}  >
   0.103 a^{-1} b^{-6} c n^{-9} > c^{\frac{1}{11}} \cdot
\frac{b^4|n|}{9.71 a} > c^{\frac{1}{11}} \,. \]
The last estimate shows that ${\lambda}_1 >0$, what we implicitly used
in (\ref{5}).

Putting these three estimates in (\ref{5}), we obtain
\[ \log z < \frac{3 \log{c} \cdot 2\log{c}}{\frac{1}{11} \log c}
=66 \log c \,.\]
Hence, $z < c^{66}$ and
\[ d=\frac{z^2 -n}{c} \leq \frac{z^2+|n|}{c}<\frac{c^{132}+c^{\frac{1}{11}}}{c}
 < c^{131}+1 .\]
\qed

Now we will develop a very useful gap principle for the
elements of a Diophantine $m$-tuple. The principle is based on the
following construction which generalizes the constructions of
Arkin, Hoggatt and Strauss \cite{A-H-S} and Jones \cite{Jon} for
the case $n=1$.

\begin{lemma} \label{l:3}
If $\{a,b,c\}$ is a Diophantine triple with the property $D(n)$ and
$ab+n=r^2$, $ac+n=s^2$, $bc+n=t^2$, then there exist integers
$e,x,y,z$ such that
\[ ae+n^2=x^2, \quad be+n^2=y^2, \quad ce+n^2=z^2 \]
and
\[ c=a+b+\frac{e}{n}+\frac{2}{n^2} (abe+rxy). \]
\end{lemma}

\pf
Define
\[ e=n(a+b+c)+2abc-2rst. \]
Then
\begin{eqnarray*}
 (ae+n^2)-(at-rs)^2 \!\!&=&\!\! an(a+b+c)+2a^2bc-2arst+n^2 \\
&&\mbox{}-a^2(bc+n)+2arst-(ab+n)(ac+n) =0.
\end{eqnarray*}
Hence we may take $x=at-rs$, and analogously $y=bs-rt$, $z=cr-st$.
We have
\begin{eqnarray*}
abe+rxy \!\!&=&\!\!
abn(a+b+c)+2a^2b^2c-2abrst \\
\mbox{}\!\!&+&\!\! abrst-a(ab+n)(bc+n)-b(ab+n)(ac+n)+rst(ab+n) \\
\!\!&=&\!\! -abcn -n^2(a+b)+rstn,
\end{eqnarray*}
and finally
\[ a+b+\frac{e}{n}+\frac{2}{n^2}(abe+rxy)=
2a+2b+c+\frac{2abc}{n}-\frac{2rst}{n}-\frac{2abc}{n}-2a-2b+\frac{2rst}{n}
=c .\]
\qed

\begin{lemma} \label{l:4}
If $\{a,b,c,d\}$ is a Diophantine quadruple with the property $D(n)$
and $|n|^3\leq a<b<c<d$, then
\[ d> \frac{3.847 \,bc}{n^2} \,.\]
\end{lemma}

\pf
We apply Lemma \ref{l:3} to the triple $\{a,c,d\}$.
Since $ce+n^2$ is a perfect square, we have that $ce+n^2\geq 0$. 
On the other hand, the assumption is that $c>|n|^3$. Hence, if $e\leq -1$, 
then $ce+n^2<-|n|^3+n^2<0$, a contradiction. Since $e$ is an integer, 
we have $e\geq 0$. If $e=0$, then $d=a+c+2s$.
If $e\geq 1$, then
\begin{equation}  \label{6}
d>a+c+\frac{2ac}{n^2}+\frac{2s\sqrt{ac}}{n^2} >\frac{2ac}{n^2} \,.
\end{equation}
(Note that if $n>0$ then $x<0$, $y<0$, and if $n<0$ and $b>|n|$ then
$x>0$, $y>0$.)

Analogously, applying Lemma \ref{l:3} to the triple $\{b,c,d\}$ we obtain
that $d=b+c+2t$ or $d>b+c+\frac{2bc}{n^2}+\frac{2t\sqrt{bc}}{n^2}$.
However, $d=b+c+2t$ is impossible since
$b+c+2t>a+c+2s$ and
\[ b+c+2t\leq b+c+2\sqrt{c(c-1)+n} <4c\leq \frac{2ac}{n^2} \,,\]
unless $a<2n^2$. But if $|n|^3\leq a <2n^2$, then $|n|=1$, $a=1$,
and in that case we have
\[ a+c+\frac{2ac}{n^2}+\frac{2s\sqrt{ac}}{n^2} >3c+2\sqrt{c(c-1)} >4c.\]

Hence we proved that
\begin{equation} \label{7}
d>b+c+\frac{2bc}{n^2}+\frac{2t\sqrt{bc}}{n^2} \,.
\end{equation}
From \cite{M-B} we know that the triples $\{1,2,3\}$ and $\{1,2,4\}$
cannot be extended to Diophantine quadruples. Thus $bc\geq 10$
and it implies
\[ t^2=bc+n\geq bc-|n|> bc-\sqrt[6]{bc} >0.853 \,bc .\]
If we put this in (\ref{7}), we obtain $d>\frac{3.847\, bc}{n^2}$.
\qed

\medskip

{\sc Proof of Theorem \ref{tm:1}.}\,\,
Assume that $\{a_1,a_2,\ldots,a_{22} \}$ has the property $D(n)$ and
$|n|^3\leq  a_1<a_2< \cdots <a_{22}$. By Lemma \ref{l:4} we find that
\begin{eqnarray*}
a_4 >\frac{a_2^2}{n^2},\quad a_5>\frac{a_2^3}{n^4}, \quad
a_6 >\frac{a_2^5}{n^8},\quad a_7>\frac{a_2^8}{n^{14}}, \\
a_8 >\frac{a_2^{13}}{n^{24}}, \quad
 a_9>\frac{a_2^{21}}{n^{40}}, \quad
a_{10} >\frac{a_2^{34}}{n^{66}},\quad a_{11}>\frac{a_2^{55}}{n^{108}}.
\end{eqnarray*}
Since $a_2>|n|^3$, we have $\frac{a_2^{55}}{n^{108}}>a_2^{11}|n|^{11}$,
and we may apply Lemma \ref{l:2} with $a=a_1$, $b=a_2$, $c=a_{11}$.
We conclude that $a_{22}\leq a_{11}^{131}$. However, Lemma \ref{l:4}
implies
\begin{eqnarray*}
a_{12} >|n| a_{11},\quad a_{13}>\frac{a_{11}^2}{|n|}, \quad
a_{14} >\frac{a_{11}^3}{n^2},\quad a_{15}>\frac{a_{11}^5}{|n|^5}, \\
a_{16} >\frac{a_{11}^8}{|n|^9},\quad a_{17}>\frac{a_{11}^{13}}{n^{16}}, \quad
a_{18} >\frac{a_{11}^{21}}{|n|^{27}},\quad
a_{19}>\frac{a_{11}^{34}}{|n|^{45}}, \\
a_{20} >\frac{a_{11}^{55}}{n^{74}},\quad a_{21}>\frac{a_{11}^{89}}{|n|^{121}},
\quad
a_{22} >\frac{a_{11}^{144}}{|n|^{197}}.
\end{eqnarray*}
Since $a_{11}>a_2^{11}|n|^{11}>n^{44}$, we obtain
\[ a_{22}>\frac{a_{11}^{144}}{|n|^{197}} \geq a_{11}^{144-\frac{197}{44}}
>a_{11}^{139}>a_{11}^{131} \,,\]
a contradiction.
\qed

\section{Small elements}

\begin{lemma}  \label{l:5}
If $\{a,b,c,d\}$ is a Diophantine quadruple with the property $D(n)$,
$|n|\neq 1$, and $n^2\leq a<b<c<d$, then $c>3.88\,a$ and $d>4.89\,c$.
\end{lemma}

\pf
We will apply Lemma \ref{l:3}. Since $b>n^2$, we have $e\geq 0$. Thus Lemma
\ref{l:3} implies that
\[ c\geq a+b+2r. \]
Since $|n|\neq 1$ we have $ab\geq 20$ and $r^2\geq ab-\sqrt[4]{ab}
>0.89 ab >0.89 a^2$. Hence, $c>3.88 \,a$.

Since $d\geq b+c+2t>a+c+2s$, from (\ref{6}) we conclude that
\[ d>a+c+\frac{2ac}{n^2}+\frac{2s\sqrt{ac}}{n^2} \,.\]
We have $ac\geq 24$ and $s^2\geq ac-\sqrt[4]{ac}> 0.9\, ac$. Therefore
\[ d>a+c+\frac{3.89\,ac}{n^2} > 4.89\,c \,.\]
\qed

{\sc Proof of Theorem \ref{tm:2}.}\,\,
We may assume that $|n|\geq 2$ since $B_1=B_{-1}=0$. Let
$\{a_1,a_2,\ldots, a_m\}$ be a Diophantine $m$-tuple with the
property $D(n)$ and
$n^2<a_1<a_2<\cdots <|n|^3$. By Lemma \ref{l:5} we have
\begin{eqnarray*}
a_3>3.88 a_1,\quad a_4>3.88\cdot 4.89a_1, \quad
\ldots \,, \quad a_m>3.88\cdot {4.89}^{m-3}a_1.
\end{eqnarray*}

Therefore
\[ 3.88\cdot {4.89}^{m-3}\cdot n^2 <|n|^3 \]
and from $m-3 <\frac{\log{\frac{|n|}{3.88}}}{\log{4.89}}$ we
obtain $m<0.65 \log{|n|} +2.24$.
\qed

\section{Very small elements}

We are left with the task to estimate the number of "very small"
elements in a Diophantine $m$-tuple. Let $\{a_1,a_2,\ldots,a_m\}$
be a Diophantine $m$-tuple with the property $D(n)$ and assume
that $a_1<a_2<\cdots <a_m\leq N$, where $N$ is a positive integer.
Let $1\leq k<m$. Then $x=a_{k+1},\,\ldots \,,\,x=a_{m}$ satisfy
the system
\begin{equation} \label{10}
a_1x+n=\Box, \quad a_2x+n=\Box, \quad \ldots \,,\quad a_kx+n=\Box,
\end{equation}
where $\Box$ denotes a square of an integer. Denote by
$Z_k(N)$ the number of solutions of system (\ref{10}) satisfying
$1\leq x\leq N$.

Motivated by the observations from the introduction of \cite{B-G-P},
we will apply a sieve method based on the following theorem of Gallagher
\cite{Gal} (see also \cite[p.29]{Hoo}):

\begin{theorem}[\cite{Gal}] \label{tm:G}
If all but $g(q)$ residue classes $({\rm mod}\,\,q)$ are removed for
each prime power $q$ in a finite set ${\cal S}$, then the number of
integers which remain in any interval of length $N$ is at most
\[ \Big( \sum_{q\in{\cal S}} \Lambda(q) -\log{N}\Big) \Big/
\Big( \sum_{q\in{\cal S}} \frac{\Lambda(q)}{g(q)} - \log{N}\Big) \]
provided the denominator is positive. Here $\Lambda(q)=\log{p}$
for $q=p^{\alpha}$.
\end{theorem}

We will use Theorem \ref{tm:G} to estimate the number $Z_k(N)$.
For this purpose, we will take
\[ {\cal S}=\{ p \,:\, \mbox{$p$ is prime, $83\leq p\leq Q$,
$\gcd(a_1a_2\cdots a_k,p)=1$} \}, \]
where $Q$ is sufficiently large. For a prime $p\in {\cal S}$ we may remove
all residue classes $({\rm mod}\,\,p)$ such that
$\Big( \frac{a_ix+n}{p} \Big)=-1$ for some $i\in \{1,\ldots,k\}$.
Here $(\frac{\cdot}{p})$ denotes the Legendre symbol.

Let $1\leq l\leq k$. Then
\begin{eqnarray*}
g(p) \!\!&\leq&\!\!
| \{ x\in {\bf F}_p \,:\, \Big( \frac{a_ix+n}{p}\Big)=0 \,\,
\mbox{or $1$, for $i=1,\ldots,l$} \} | \\
\!\!&\leq&\!\!
l+ | \{ x\in {\bf F}_p \,:\, \Big( \frac{x+n\overline{a_i}}{p}\Big)=
\Big( \frac{\overline{a_i}}{p}\Big), \,\,
\mbox{for $i=1,\ldots,l$} \} | \,.
\end{eqnarray*}
Here $a_i\overline{a_i}\equiv 1 \!\!\pmod{p}$. Using estimates 
for character sums (see \cite[p.325]{L-N}), we obtain
\[ g(p)\leq l+\frac{p}{2^l} +\Big( \frac{l-2}{l}+\frac{1}{2^l}\Big)\sqrt{p}
+\frac{l}{2} \,.\]

Assume that $k=\lfloor \log_{2}{Q} \rfloor$. We may take
$l= \lfloor \log_{2}{p} \rfloor$. Then we have
\[ \frac{p}{2^l}+\frac{\sqrt{p}}{2^l} +\frac{3l}{2}<2+\frac{2}{\sqrt{p}}+
\frac{3\log_{2}{p}}{2} <\sqrt{p} \]
for $p\geq 179$. Hence
\[ l+\frac{p}{2^l}+\Big( \frac{l-2}{2}+\frac{1}{2^l}\Big)\sqrt{p}
+\frac{l}{2}<\frac{l}{2}\sqrt{p}< \frac{\log_{2}{p}}{2}\sqrt{p}
<0.722 \,\sqrt{p}\,\log{p} \]
for $p\geq 179$, and we may check directly that
$l+\frac{p}{2^l}+\Big(\frac{l-2}{2}+\frac{1}{2^l}\Big) +\frac{l}{2}
< 0.722 \,\sqrt{p}\log{p}$ \,for $83\leq p\leq 173$.
Therefore we proved that
\[ g(p)<0.722 \,\sqrt{p} \log{p}. \]

\medskip

By Theorem \ref{tm:G}, we have $Z_k(N)\leq \frac{E}{F}$, where
\[ E=\sum_{p\in {\cal S}} \log{p} -\log{N}, \quad
F=\sum_{p\in {\cal S}} \frac{1}{0.722\,\sqrt{p}} -\log{N}. \]
By \cite[Theorem 9]{R-S}, we have $E<\sum_{83\leq p\leq Q} \log{p}
<\theta(Q) < 1.01624\,Q$.

Assume that at least $\frac{4}{5} \pi(Q)$ primes less than $Q$
satisfy the condition $\gcd(a_1a_2\cdots a_k,p)=1$. Then we have
\begin{eqnarray} \label{11}
F \!\!&\geq&\!\! \frac{1}{0.722\,\sqrt{Q}} |{\cal S}| -\log{N}
\geq \frac{1}{0.722\,\sqrt{Q}} \Big( \frac{4}{5}\pi(Q) -23\Big)-\log{N}
\nonumber \\
\!\!&>&\!\! 1.108 \frac{\sqrt{Q}}{\log{Q}}-\frac{31.86}{\sqrt{Q}}-\log{N} \,.
\end{eqnarray}
Since $F$ must be positive in the applications of Theorem \ref{tm:G},
we will choose $Q$ of the following form
\begin{equation} \label{12}
Q= c_1 \cdot \log^2{N} \cdot (\log\log{N})^2,
\end{equation}
where $c_1$ is a constant.

We have to check whether our assumption is correct. Suppose that
$a=a_1a_2\cdots a_k$ is divisible by at least one fifth of the primes
$\leq Q$. Then $a\geq p_1p_2\cdots p_{\lceil \frac{1}{5} \pi(Q) \rceil}$,
where $p_i$ denotes the $i^{\rm th}$ prime. By \cite[p.69]{R-S},
we have
\[ p_{\lceil \frac{1}{5} \pi(Q) \rceil} > \frac{1}{5}\pi(Q) \log(\frac{1}{5}
\pi(Q))  >\frac{1}{5}\frac{Q}{\log{Q}} \log\Big(\frac{1}{5}\frac{Q}{\log{Q}}
\Big) := R.\]
Therefore, by \cite[p.70]{R-S},
\[ \log{a}>\sum_{p\leq R} \log{p} >R\Big(1-\frac{1}{\log{R}} \Big). \]
Assume that $Q\geq 2\cdot {10}^4$. Then $\frac{1}{5}\frac{Q}{\log{Q}}
>Q^{0.605}$ and $R>0.128\,Q$. Furthermore,
$\log{R}>7.793$ and therefore
\[ \log{a}>0.105\,Q.\]
On the other hand, $a<N^k$ and $\log{a}<k\log{N}\leq \log_{2}{Q}\log{N}$.

Assume that $N\geq 1.6\cdot {10}^5$ and $c_1\leq 80$. Then
$Q\leq \log^{4.498}{N}$. In order to obtain a contradiction,
it suffices to check that
\[ 0.105\,c_1\, \log^{2}{N} \,(\log\log{N})^2 >
\frac{4.498}{\log{2}} \log{N} \cdot \log\log{N} \]
or
\[ c_1\,\log{N} \,\log\log{N} >61.81, \]
and this is certainly true for $N\geq 1.6\cdot {10}^5$ if we
choose $c_1\geq 2.08$.

\medskip

Thus we may continue with estimating the quantity $F$. We are working
under assumptions that (\ref{12}) holds with $2.08\leq c_1\leq 80$,
$Q\geq 2\cdot{10}^4$ and $N\geq 1.6\cdot{10}^5$. We would like to have
the estimate of the form
\begin{equation} \label{13}
F>\frac{\sqrt{Q}}{c_2\,\log{Q}} \,.
\end{equation}
This estimate will lead to
\begin{equation} \label{14}
Z_k(N) < 1.01624\, c_2 \sqrt{Q}\log{Q} < 4.572\, c_2\sqrt{c_1} \log{N}
(\log\log{N})^2.
\end{equation}

In order to fulfill (\ref{13}), it suffice to check
\[ \frac{31.86}{\sqrt{Q}}+\log{N} <\frac{\sqrt{Q}}{\log{Q}}
\Big( 1.108-\frac{1}{c_2}\Big).  \]
Since $Q>2\cdot {10}^4$ we have $\frac{31.86}{\sqrt{Q}} <0.016\,\frac{\sqrt{Q}}
{\log{Q}}$. Furthermore,
\[ \frac{\log{N}\log{Q}}{\sqrt{Q}} < \frac{4.498\,\log{N}\log\log{N}}
{\sqrt{c_1}\log{N}\log\log{N}} =\frac{4.498}{\sqrt{c_1}} \,.\]
Hence $c_2>1 \Big/ (1.092 -\frac{4.498}{\sqrt{c_1}})$. Thus if we choose
$c_1=68$, then we may take $c_2=1.83$ and from (\ref{14}) we obtain
\begin{equation} \label{15}
Z_k(N)<69\,\log{N}(\log\log{N})^2.
\end{equation}
Note that with this choice of $c_1$, $N\geq 1.6\cdot{10}^5$ implies
$Q>60222>2\cdot{10}^4$.

\medskip

{\sc Proof of Theorem \ref{tm:3}.}\,\,
Let $\{a_1,a_2,\ldots,a_m\}$ be a Diophantine $m$-tuple with the
property $D(n)$ and $a_1<a_2<\cdots <a_m\leq n^2$. Then for any
$k\in \{1,2,\ldots,m\}$ we have
\[ m\leq k+ Z_k(n^2). \]
Let $k=\lfloor \log_2{Q} \rfloor$, where $Q=68\,\log^2{n^2}\,
(\log\log{n^2})^2$. Since $|n|\geq 400$, we have $n^2\geq 1.6\cdot{10}^5$ and
we may apply formula (\ref{15}) to obtain
\begin{equation} \label{16}
Z_k(n^2) <69\,\log{n^2} (\log\log{n^2})^2 <
265.55\,\log|n| \,(\log\log|n|)^2,
\end{equation}
Furthermore,
\begin{equation} \label{17}
k<\frac{1}{\log_2}\,\log(\log^{4.489}{n^2}) < 9.01 \,\log\log|n| ,
\end{equation}
and combining (\ref{16}) and (\ref{17}) we finally obtain
\[ m < 265.55 \,\log{|n|}(\log\log{|n|})^2 +9.01 \,\log\log{|n|}. \]
\qed

\medskip

\begin{remark} \label{r:k}
{\rm In \cite{Kat} Katalin Gyarmati recently considered the more 
general problem. She estimated $\min \{|{\cal A}|, |{\cal B}|\}$, where 
${\cal A}, {\cal B} \subseteq \{1,2,\ldots,N\}$ satisfy the condition 
that $ab+n$ is a perfect square for all $a\in {\cal A}$, $b\in {\cal B}$. 
Using her approach, it can be deduced that if $\{a_1,a_2,\ldots,a_m\}$ 
has the property $D(n)$, where $n>0$ and $a_1<a_2<\cdots <a_m\leq N$, 
then $m\leq 2n \log{N}$. This yields $C_n\leq 4n\log{n}$ for $n\geq 2$. }  
\end{remark}

\begin{remark} \label{r:rss}
{\rm Let us mention that Rivat, S\'ark\"ozy and Stewart \cite{RSS} 
recently used Gallagher's "larger sieve" method is estimating the size 
of a set $Z$ of integers such that $z+z'$ is a perfect square whenever 
$z$ and $z'$ are distinct elements of $Z$. They proved that if 
$Z\subset \{1,2,\ldots, N\}$, where $N$ is greater that an effectively 
computable constant, then $|Z|< 37 \,\log N$. 

Largest known set with the above property is a set with six elements 
found by J.~Lagrange \cite{Lag}. Maybe this may be compared with our 
situation where the largest known Diophantine $m$-tuples are Diophantine 
sextuples found by Gibbs \cite{Gibbs1,Gibbs2}. } 
\end{remark}

\medskip

{\sc Proof of Theorem \ref{tm:4}.}\,\,
Since $M_n\leq A_n+B_n+C_n$, the second part of the theorem follows
directly from Theorems \ref{tm:1}, \ref{tm:2} and \ref{tm:3}.

For $|n|\leq 400$, Theorem \ref{tm:2} gives $B_n\leq 6$.
It is easy to verify with a computer that for $|n|\leq 400$
it holds $C_n\leq 5$. More precisely, $C_n=5$ if and only if
$n\in \{-299,-255,256,400\}$. These two estimates together with
Theorem \ref{tm:1} imply $M_n\leq 32$.
\qed

\section{Concluding remarks}
It is not surprising that in Theorem \ref{tm:4} the main contribution
comes from $C_n$. Namely, if we define $C=\sup \{C_n \,:\, n\in
\Z\setminus \{0\} \}$, then we have $M=C$. Indeed, if
$\{a_1,a_2,\ldots,a_m\}$ is a Diophantine $m$-tuple with the property
$D(n)$, then $\{a_1c,a_2c,\ldots,a_mc\}$ has the property $D(nc^2)$
and for sufficiently large $c$ we have $a_ic\leq (nc^2)^2$, $i=1,2,\ldots,m$.
It means that in order to prove $M<\infty$, it suffices to prove
$C<\infty$. The above argumentation shows that it suffice to prove
that for some $\varepsilon>0$ it holds
\[ \sup_{n\neq 0}\,\sup \{|S\cap [1,n^{0.5+\varepsilon}]| \,:\,
\mbox{$S$ has the property $D(n)$} \}<\infty. \]

We may define also $A=\sup \{A_n \,:\, n\in \Z\setminus \{0\} \}$ and
$B=\sup \{B_n \,:\, n\in \Z\setminus \{0\} \}$. Gibbs' example mentioned
in introduction shows that $C\geq 6$ and $M\geq 6$.
If $n=a^2$, $a\geq 5$, then $B_n\geq 3$ since $\{a^2+1,\,a^2+2a+1,\,
4a^2+4a+4\}$ has the property $D(a^2)$. Hence $B\geq 3$.
Finally, since $\{k,\,k+2,\,4k+4,\,16k^3+48k^2+44k+12\}$ has the property
$D(1)$ we have $A\geq A_1\geq 4$.

{\footnotesize
{\sc Department of Mathematics, University of Zagreb, Bijeni\v cka cesta 30,
10000 Zagreb, Croatia}

{\em E-mail address}: {\tt duje@math.hr} }

\end{document}